\input amstex
\documentstyle{amsppt}

\define\cal{\Cal}
\define\C{\Bbb C}
\define\a{\alpha}

\define\De{\Delta}

\define\Z{\Bbb Z}
\define\ra{\rightarrow}
\define\N{\Bbb N}
\define\Hom{\hom}

\topmatter
\title Quantum groups and representations with highest weight 
\endtitle

\author Joseph Bernstein and Tanya Khovanova \endauthor

\address Dept. Math. Tel-Aviv University, Ramat Aviv, 69978 Israel
\endaddress
\email{bernstei\@math.tau.ac.il}\endemail

\address Dept. Math. Princeton University, Princeton, NJ 08544
\endaddress
\email{tanya\@math.princeton.edu}\endemail

\abstract
We consider a special category of Hopf algebras, depending on parameters
$\Sigma$ which possess
properties similar to the category of representations of simple Lie group
with highest weight $\lambda$. We connect quantum groups to minimal objects
in this categories---they correspond to irreducible representations
in the category of representations with highest weight $\lambda$. Moreover,
we want to correspond quantum groups only to finite dimensional irreducible
representations. This gives us a condition for $\lambda$: $\lambda$---
is dominant means the minimal object in the category of representations
with highest weight $\lambda$ is finite dimensional. We put
similar condition for $\Sigma$. We call $\Sigma$ dominant if
the minimal object in corresponding category has polynomial growth.
Now we propose to define quantum groups starting from dominant
parameters $\Sigma$.
\endabstract

\date April 8, 1997 
\enddate

\endtopmatter

\document
\head
{1. Definitions and examples}
\endhead
\subhead{1.1 Torus}\endsubhead
Let us fix an $n$-dimensional torus $H$ 
(i.e. an algebraic group isomorphic to $\C^{*n}$). We denote
by $S$ the Hopf algebra of regular functions on $H$.
Let $\Lambda$ be the lattice of characters of $H$. Then 
$\Lambda \subset S$ is a basis in $S$. The dual algebra $S^*$ 
can be realized as the algebra of functions on the lattice
$\Lambda$.
We denote by $\hat H$ the group algebra of $H$: 
for
any element $h \in H$ we denote the corresponding generator
in $\hat H$ by $\hat h$ ($\hat h_1 \hat h_2 = \widehat {h_1h_2}$). 
The Hopf algebra $\hat H$ is a subalgebra in
$S^* \ (\hat h(\lambda) = \lambda(h))$ on which the comultiplication 
is well defined:
$\De \hat h = \hat h \otimes \hat h$. 

The Hopf algebra $\hat H$ is
too small for some of our future purposes. In order to be able to
define the comultiplication on $S^*$ we have to complete $S^* \otimes
S^*$ to $(S \otimes S)^*$. We need the comultiplication to define an
action of $S^*$ on the tensor product of two $S^*$-modules.
We are interested only in those $S^*$-modules $W$ for which 
the $S^*$-module
structure is inherited from some $S$-comodule structure.
That means, we should consider only $S^*$-modules $W$ which are
algebraic representations of $H$. Under this condition the
completed comultiplication on $S^*$ would allow us to define
the action of $S^*$ on the tensor product of two such $S^*$-modules 
(see \cite{B-Kh}).

We set $S^* \hat \otimes S^* :=
(S \otimes S)^*$. The algebra $S^* \hat \otimes S^*$ could be realized
as the algebra of all functions on the lattice $\Lambda \oplus
\Lambda$. If $ f \in S^*$ then $\De f(\lambda _1,\lambda_2) =
f(\lambda _1 + \lambda_2)$.
\medskip
For convenience, let us denote by $S^\star$ either $\hat H$, or
$S^*$ with restrictions described above 
(or, any other suitable representative of a dual Hopf algebra of $S$)

\subhead{1.2 Datum}\endsubhead
Given $H$ our datum 
$\Sigma$ would be two 
finite sets of size $m$: $\{\a_1, \a_2,\allowmathbreak ...,\a_m\} 
\allowmathbreak 
\subset 
\Lambda \setminus \{0\}$ ---
the set of non-zero characters; and $\{\gamma_1, \gamma_2,...,
\gamma_m\} \subset H$ --- 
the set of points of the torus. We denote a generator
$\hat \gamma_k \in \hat H$ corresponding to the point $\gamma_k \in H$
by $K_k$. We denote 
$\a_i(\gamma_j)$ by $q_{ij}$.

\subhead{1.3 Tetramodules}\endsubhead
Using our datum we can construct 
an $S$-tetra\-module $T$ and an $\hat H$-tetramodule $V$ (we can consider
$V$ as an $S^*$-tetramodule, see above). For
definition of tetramodule (Hopf bimodule etc.) (see \cite{B-Kh} and
references there). Informally,
S-tetramodule is an $S$-bimodule and $S$-bicomodule with some
natural axioms. The tetramodule $T$ is generated over $S$ by its
space of right $H$-invariants. The elements $t_i \ \ \ 1 \le i \le m$ would
generate a linear basis in this space. Then we describe an $S$-tetramodule
structure of $T$ as follows:
$$\De t_i = t_i \otimes 1 + \a_i \otimes t_i$$
$$s t_i s^{-1} = s(\gamma_i)t_i \ \ \ \text{for}\ \  s \in \Lambda.$$

Analogously, an $\hat H$-tetramodule $V$ is generated over $\hat H$
by elements $E_i \ \ 1 \le i \le m$ with the following tetramodule
structure:
$$\De E_i = E_i \otimes 1 + K_i \otimes E_i$$
$$\hat h E_i \hat h^{-1} = \a_i(h)E_i \ \ \text{for} \ \ \hat h 
\in \hat H \subset S^*.$$

\subhead{1.4 Categories}\endsubhead
Given an $S$-tetramodule $T$
denote by ${\cal H} (S,T)$ the category of $\Z_+$-graded
Hopf algebras $B$ such that $B_0 = S, \ B_1 = T$; and $B$ supplies
$T$ with the given $S$-tetramodule structure.

Given our datum we have two categories---${\cal H} (S,T)$ and
${\cal H}(S^\star,V)$. By latter category we mean either
${\cal H}(\hat H,V)$, or ${\cal H} (S^*,V)$ with restrictions
discussed above.

\subhead{1.5 Examples}\endsubhead
1. Let $G$ be a simple Lie group, and $H$---its
Cartan subgroup. Consider datum $\Sigma$ depending on a parameter
$q$. We put $\Sigma$ to be
the set $\{\a_i\}$ of simple roots and
the set $\{\gamma_i\}$ such that $\a_i(\gamma_j) = q^{<\a_i,\a_j>}$,
where the elements $<\a_i,\a_j>$ form a
Cartan matrix of $G$. Then the universal enveloping algebra of a
Borel subalgebra $B_+$ of a quantum group $G_q$ is an object in
${\cal H}(S^*,V)$. Actually, when $q \ne 1$ the Hopf algebra
$U(B_+)$ is an object in ${\cal H}(\hat H,V)$, but the limit
$q \ra 1$ should be considered in the bigger algebra.

2. Let $G$ be a reductive algebraic group; $H$---its Cartan subgroup.
Let $A = \C [G]$ be the Hopf algebra of regular functions on $G$ and
$I$ the Hopf ideal of functions equal to $0$ on $H$. Then $S = \C[H]$
equals $A/I$; and the adjoint graded Hopf algebra $gr A$ (with respect to $I$)
is the object in ${\cal H} (S,T)$, where $T$ is described by datum
$\{\a_i\}$---the set of non-zero roots and all $\gamma_i$ equals $1$.

3. Let $G_q$ be a quantum deformation of a simple Lie group $G$.
By this we mean a flat family of Hopf algebras $A_q$ which are deformations
of $A = \C[G]$. It can be shown that we can flatly deform an
ideal $I$ (see Example 2) so that the family of quotient Hopf algebras
$H_q = A_q/I_q$ is constant and equals $S = \C [H]$. It is
easy to see that the adjoint graded algebra $gr A_q$ for generic $q$
is the object in ${\cal H} (S,T)$, where $T$ is defined by following
datum $\Sigma$: the set of characters is the set $\{\a_i, -\a_i\}$, where
$\a_i$ are all simple roots, the set of points is $\{\gamma_i, \gamma_i\}$
such that $\a_i(\gamma_j)$ are defined by Cartan matrix of $G$.
\medskip
Given a quantum group $G_q$ we can construct
a Hopf algebra $B_q^1 \in {\cal H}(S^\star,V)$ for some datum $\Sigma^1(G_q)$,
see example 1. Also, there is another construction (from example 3)
of a Hopf algebra $B^2_q \in {\cal H}(S,T)$ for some other datum
$\Sigma^2(G_q)$. We suggest that we have some construction (of example 1
or 3, or maybe similar) such that we can describe any quantum group
$G_q$ in terms of some Hopf algebra $B_q$, which is an object in
the category ${\cal H}(S,T)$ or ${\cal H}(\hat H,V)$ (resp. 
${\cal H}(S^*,V)$) for some
datum $\Sigma$. 

\subhead{1.6 Our goals}\endsubhead
Given a quantum group $G_q$ we constructed a $\Z_+$-graded Hopf algebra
in some category ${\cal H}(S,T)$ for some datum $\Sigma$ (see 1.5). 
We would like to answer the following questions:

1) Given datum $\Sigma$ how we can distinguish an object in the
category ${\cal H}(S,T)$ which could correspond to a quantum group.

2) What properties should $\Sigma$ satisfy in order to supply
the category ${\cal H}(S,T)$ with an object which correspond
to some quantum group.

\head{2. The category of modules with highest weight}
\endhead

\subhead{2.1}\endsubhead
Our intuition and constructions are partly based on the idea
from our previous paper \cite{B-Kh} of deep parallelism of properties of
the category ${\cal H}(S,T)$ (resp. ${\cal H}(S^*,V))$ and the category of
modules with highest weight $\lambda$.

Our datum now are the simple Lie group $G$ and the
weight $\lambda , \ \lambda \in \frak H^*$. We consider the
category ${\cal O}^\prime$
of modules over $G$ from the category ${\cal O}$, such that all their
weights belong to the set $\{\lambda-\Gamma_+\}$, where 
$\Gamma_+$ is generated over $\N$ by simple roots.

\subhead{2.2}\endsubhead
Consider the functor $J$: ${\cal O}^\prime \ra {\cal V}$, where
we denote the category of vector spaces by ${\cal V}$,
such that to any module $M$ we correspond the vector
space $X$ of vectors of weight $\lambda$.
\medskip
\proclaim{Lemma} The functor $J$ possess the left adjoint
functor $F$: ${\cal V} \to {\cal O}^\prime$. \endproclaim
\medskip
Consequently, for any $X \in {\cal V}$ we can construct 
a module $FX$, such that
$$\Hom_{{\cal O}^\prime}(FX,M)=\Hom_{\cal V}(X,JM).\leqno(1)$$
\medskip
\example{Example} If $X=\C$ then $FX=M_\lambda$---the Verma module
with highest weight $\lambda$ which possess the fixed vector
of weight $\lambda$. The equality (1) means that for any
module $M$ and for any morphism
$\C \ra JM$ 
we can construct a unique morphism
$M_\lambda \ra M$ which is identical on the image of $\C$.
Denote by ${}^\C{\cal M}$ the category of modules from ${\cal O}^\prime$
together with the fixed morphism $\C \to JM$.
Then 
we can rephrase our example, saying
that the Verma module is the initial object
in the category ${}^\C{\cal M}$ (compare \cite{B-Kh}).
\endexample
\medskip
\demo{Proof} The example above gives us a proof when $\dim X=1$. 
The proof for any $X$ could be easily modified from this
example.
\enddemo
\medskip
\proclaim{Lemma} The functor $J$ possess the right adjoint
functor $H$: ${\cal V} \to {\cal O}^\prime$.
\endproclaim
\medskip
Consequently, for any $X$ in ${\cal V}$ we can construct a
module $FX$ such that:
$$\Hom_{\cal V}(JM, X)=\Hom _{{\cal O}^\prime}(M, HX).\leqno(2)$$
\medskip \example{Example} If $X$=$\C$ then $FX=\delta_\lambda$---the
contragredient Verma module with fixed highest covector of weight
$\lambda$. The equality above means that for any module $M$ and
any morphism $JM \to \C$ 
we can construct a unique morphism $M \to \delta_\lambda$
which is identical on the fixed covector. 
Denote by ${\cal M}^\C$ the category of modules from ${\cal O}^\prime$
together with the fixed morphism $JM \to \C$.
Then we can rephrase our example, saying that the contragredient
Verma module $\delta_\lambda$ is the final object in the
category ${\cal M}^\C$ (compare \cite{B-Kh}).
\endexample
\medskip
\demo{Proof}In the category ${\cal O}^\prime$ there is a
natural duality: $M \to M^\star$, where we define each
weight subspace of $M^\star$ as a regular dual to corresponding
weight subspace in $M$: $(M^\star)_\eta=(M_\eta)^*$.
The action of the group $G$ is defined naturally. After this
remark the existance of final object in ${\cal M}^\C$ 
becomes automatic and the construction of the contragredient Verma
module becomes trivial: $\delta_\lambda=(M_\lambda)^\star$.
\enddemo

\subhead{2.3 Shapovalov map}\endsubhead
Consider the category ${\cal M}$ of modules $M \in{\cal O}^\prime$
such that $\dim (JM)=1$ together with fixed isomorphism
between $JM$ and $\C$. This category is the subcategory
in ${}^\C{\cal M}$ and in ${\cal M}^\C$. As $M_\lambda \in {\cal M}
\subset {}^\C{\cal M}$
it is an initial object in ${\cal M}$, analogously, $\delta_\lambda$
is a final object in ${\cal M}$. Therefore, there exists a canonical
map $Sh$: $M_\lambda \to \delta_\lambda$ which is called a Shapovalov
map. (As $\delta_\lambda \subset M_\lambda^*$ the Shapovalov map
defines the Shapovalov form on $M_\lambda$).
\medskip
Going back to functors, if we put $M=FX$ into $(1)$, we get
$$\Hom_{{\cal O}^\prime}(FX,FX)=\Hom_{\cal V}(X,JFX).$$
The identity isomorphism on the left corresponds to 
a canonical element $j$ on the right which is called
the adjunction map: $j$: $X \to JFX$.
It is easy to check that in the category ${\cal M}$
the adjunction map $j$: $X \to JFX$ is an isomorphism $\forall X
\in {\cal V}$. If we put instead of $M$ the module $FX$ in $(2)$ we
would get:
$$\Hom_{{\cal O}^\prime}(FX,HX)=\Hom_{\cal V}(JFX,X).$$
For the category ${\cal M}$
the identity element in the right hand side (which is inverse
of the adjunction map $j$) would correspond
to a canonical morfisms of functors on the left:
\proclaim{Lemma} In the category ${\cal M}$ there 
exists a canonical morphism of
functors $F \to H$.
\endproclaim
\medskip
We denote $L_\lambda$ the image of the Shapovalov map in
$\delta_\lambda$. The module
$L_\lambda$ is an irreducible representation of $G$ with the highest 
weight $\lambda$. The module $L_\lambda$ is the minimal object
in ${\cal M}$. That means that for any module $B \in {\cal M}$
there exists a submodule $B^\prime \subset B$ and a canonical
epimorphism $B^\prime \to L_\lambda$.

\subhead{2.4 Point}\endsubhead
We would keep in mind for the next
section that minimal object in the category ${\cal M}$
is of importance. Also, if we consider the category ${\cal M}$
as a function of $\lambda$ we may say that we are interested in
those $\lambda$'s when the dimension of minimal object in
${\cal M}$ is dropping significantly. (Weight $\lambda$ is
called dominant if $L_\lambda$ is finite dimensional).

\head{3. Our categories}\endhead
\subhead{3.1 Parallel construction}\endsubhead
Let us fix $S$---the Hopf algebra of functions over torus.
We denote by ${\cal T}$ the category of tetramodules over $S$.
We denote by ${\cal A}$ the category of $\Z_+$-graded Hopf
algebras with $0$-graded subalgebra isomorphic to $S$ and
the subspace of each grade is finitely generated over $S$. Consider
the functor $J$: ${\cal A} \to {\cal T}$ which corresponds
to a given Hopf algebra $A$ a $1$-graded subspace of $A$ with
inherited $S$-tetramodule structure \cite{B-Kh}.
\proclaim{Lemma} The functor $J$ possess the left adjoint functor
$F$: ${\cal T} \to {\cal A}$.
\endproclaim
\medskip
\demo{Proof} 
Given an $S$-tetramodule $T$ we can construct 
a $\Z_+$-graded algebra $A$ such that $A_0 = S,
\ A_1 = T$, and the algebra $A$ is freely
generated as an algebra by $S$ and $T$ and $A$ supplies
$T$ with the given $S$-bimodule structure. 
We can construct a comultiplication on $A$
by multiplicativity. Being universal as an algebra,
the Hopf algebra $A$ remains universal as a Hopf algebra.
\enddemo
\medskip
\proclaim{Lemma} The functor $J$ possess the right adjoint functor
$F$: ${\cal T} \to {\cal A}$.
\endproclaim
\medskip
\demo{Proof} 
Given an $S$-tetramodule $T$ we can construct 
a $\Z_+$-graded coalgebra $A$ such that $A_0 = S,
\ A_1 = T$, and the coalgebra $A$ is freely
generated as a coalgebra by $S$ and $T$ and $A$ supplies
$T$ with the given $S$bicomodule structure. 
We can construct a multiplication on $A$
by comultiplicativity. Being universal as a coalgebra,
the Hopf algebra $A$ remains universal, when the Hopf algebra
structure is added \cite{B-Kh}. 
\enddemo
\medskip
\demo{Remark} We can construct a category ${\cal A}^\star$ of $\Z_+$-graded
Hopf algebras dual to the category ${\cal A}$. To each Hopf algebra
$A \in {\cal A}$ we correspond a Hopf algebra $A^\star \in {\cal
A}^\star$ such that $A^\star$ as an algebra is a subalgebra in $A^*$
and each component $(A^\star)_n$ of $A^\star$ is isomorphic
as a vector space to $S^\star \otimes M^*$, when $A_n = S\otimes M$.
After that our universal coalgebra in ${\cal A}$ corresponds
to universal algebra in ${\cal A}^\star$.
\enddemo
\medskip
As in section 2, let us consider the category ${\cal H}(S,T)$ 
and for each element $A \in {\cal H}$ let us fix an
isomorphism of $A/A_2$ with $S\oplus T$.
It is easy to see that the category ${\cal H}(S,T)$ is  similar
to the category ${\cal M}$ in chapter $2$.

In particular, there exists 
an initial object in the category ${\cal H}(S,T)$. It is called
a universal algebra. We denote it by $B^i(S,T)$. It corresponds
to Verma module $M_\lambda$.

There exists a final object in the category ${\cal H}(S,T)$.
It is called a universal coalgebra. We denote it by $B^f(S,T)$.
It corresponds to contragredient module $\delta_\lambda$.

Hence, there is a canonical map $Sh: B^i \ra B^f$ which is
analogue of the Shapovalov map. We would denote the image of 
this map by $B^m(S,T)$. This would be the minimal object
in the category ${\cal H}(S,T)$. The Hopf algebra $B^m$ corresponds
to an irreducible representation in our parallelism.

\subhead{3.2 Answers to questions}\endsubhead
Constructing a parallelism between
the categories ${\cal M}$ and ${\cal H}(S,T)$ we now can
use our intuition in ${\cal M}$ to understand what is important in
${\cal H}(S,T)$. Now we are ready to answer to our first
question. If quantum groups corresponds to datum $\Sigma$, then
it corresponds to a minimal object $B^m$ in the category
${\cal H}(S,T)$, where the category ${\cal H}(S,T)$ is constructed
by our datum.
\medskip
We would call our datum $\Sigma$ dominant if the minimal object
$B^m$ in the category ${\cal H}(S,T)$ has polynomial growth over $S$.
\medskip
As an answer to our second question,
we suggest that quantum groups correspond to Hopf algebras which
are minimal objects constructed from dominant datum.

\enddocument